\documentclass[11pt]{amsart}

\usepackage{cancel}

\newcommand{\sout}[1]{ $ \cancel{#1} $ }

\newcommand{\floor}[1]{\lfloor{#1}\rfloor}
\newcommand{\Sprimes}{\mathcal{S}}
\newcommand{\Wprimes}{\mathcal{W}}

\newtheorem{theorem}{Theorem}
\newtheorem{conjecture}{Conjecture}
\newtheorem{lemma}{Lemma}

\newcommand{\ktuple}{ (f_1(x), \ldots, f_k(x) )}

\begin{document}
\bibliographystyle{plain}

\title{Two Algorithms to Find Primes in Patterns}

\author{Jonathan P.~Sorenson}
\address{Computer Science and Software Engineering, 
  Butler University, Indianapolis, IN 46208 USA}
\email{sorenson@butler.edu}

\author{Jonathan Webster}
\address{Mathematics, Statistics, and Actuarial Science,
  Butler University, Indianapolis, IN 46208 USA}
\email{jewebste@butler.edu}

\date{\today}

\begin{abstract}
  Let $k\ge 1$ be an integer, and
    let $P=\ktuple$ be $k$
    admissible linear polynomials over the integers,
    or \textit{the pattern}.  
  We present two algorithms that find all integers $x$ where
    $\max{ \{f_i(x) \} } \le n$ and all the $f_i(x)$ are prime.
\begin{itemize}
\item
  Our first algorithm takes 
    $O_P(n/(\log\log n)^k)$ arithmetic operations using $O(k\sqrt{n})$ space.
\item
  Our second algorithm takes slightly more time,
    $O_P(n/(\log \log n)^{k-1})$ arithmetic operations,
    but uses only $n^{1/c}$ space, for $c>2$ a fixed constant.
  The correctness of this algorithm is unconditional, but
    our analysis of its running time depends on 
    two reasonable but unproven conjectures.
\end{itemize}
We are unaware of any previous complexity results for this problem
  beyond the use of a prime sieve.

We also implemented several parallel versions of our second algorithm
  to show it is viable in practice.
In particular, we found some new Cunningham chains of length 15,
  and we found all quadruplet primes up to $10^{17}$.

\end{abstract}


\maketitle

\section{Introduction}

Mathematicians have long been interested in prime numbers
  and how they appear in patterns.
(See, for example, \cite[ch.~A]{UPINT}.)
In this paper, we are interested in the complexity of
  the following algorithmic problem:
\begin{quote}
  Given a pattern and a bound $n$, find all primes $\le n$ that fit the
  pattern.
\end{quote}
To address this, first we will discuss and define a pattern of primes,
  then we will look at what is known about the
   distribution of primes in patterns to see what we can
   reasonably expect for the complexity of this problem,
and finally we will discuss previous work and state our new results.

\subsection{Prime Patterns}

Perhaps the simplest of patterns of primes are the twin primes,
  which satisfy the pattern $(x,x+2)$ where both $x$ and $x+2$ are prime.
Examples include 59,61 and 101,103.

We can, of course, generalize this to larger patterns.
For example, prime quadruplets have the form
  $(x, x+2, x+6, x+8)$,
and examples include 11,13,17,19 and 1481,1483,1487,1489.

Larger patterns of primes of this type are
  called \textit{prime $k$-tuples}.
If the $k$-tuple has the smallest possible
  difference between its first and last primes (its \textit{diameter}),
  it is also called a \textit{prime constellation}.
The pattern 
  $(x, x+2, x+6, x+8)$ is a constellation, as $8$ is the smallest
  possible diameter for a pattern of length 4.
There are two constellations of length 3:
  $(x,x+2,x+6)$ and $(x,x+4,x+6)$.
See, for example, \cite[\S1.2.2]{CP} or \cite[ch.\ 3]{Riesel}.

Sophie Germain studied the pattern $(x, 2x+1 )$, which was later
  generalized to \textit{Cunningham chains} of two kinds.
Chains of the \textit{first kind}
  have the pattern $(x, 2x+1, 4x+3, 8x+7, \ldots)$,
and chains of the \textit{second kind}
  have the pattern $(x, 2x-1, 4x-3, 8x-7, \ldots)$.

Chernick \cite{Chernick39} showed that any prime pattern of the form
  $(6x+1, 12x+1, 18x+1)$ gives
  a Carmichael number composed of the product of these three primes.

Let $k>0$ be an integer.
A \textit{prime pattern} of size $k$ is a list of $k$ 
  linear polynomials over the integers with positive leading coefficients,
  $\ktuple$.
A pattern of size $k$ is \textit{admissible} if 
  for every prime $p\le k$, there is an integer $x$ such that
  $p$ does not divide any of the $f_i(x)$.
\nocite{GR98}
For an algorithm to test for pattern admissibility,
  see \cite[pp. 62--63]{Riesel}.

We restrict our notion of pattern to linear polynomials in this paper.

\subsection{The Distribution of Primes in Patterns}

The unproven twin prime conjecture 
  states there are infinitely many twin primes.
Yitang Zhang \cite{Zhang14} recently showed that there is a 
  positive integer $h$
  such that the pattern $(x,x+h)$ is satisfied by infinitely many primes.

The Hardy-Littlewood $k$-tuple conjecture \cite{HL23} implies that
  each pattern, with leading coefficients of 1, 
  that is \textit{admissible}, will be satisfied
  by primes infinitely often.
Further, the conjecture implies that the number of primes $\le n$ 
  in such a pattern of length $k$ is roughly proportional to $n/(\log n)^k$.

For twin primes, then, the Hardy-Littlewood conjecture gives an estimate of
  $$ 2C_2 \frac{n}{(\log n)^2} $$
for the number of twin primes $\le n$, where $C_2\approx 0.6601\ldots$
  is the twin primes constant.
Brun's theorem gives an $O(n/(\log n)^2)$ upper bound for the
  number of twin primes $\le n$.
For every $k$-tuple with $k\ge 2$, there is a 
  corresponding conjectured constant of proportionality
  that depends on the pattern and on $k$,
  generalizing the twin primes constant.

Dickson's conjecture \cite{Dickson04} states that there are infinitely many
  primes satisfying any fixed admissible pattern,
  even with leading coefficients $> 1$.
Thus, it applies to fixed-length Cunningham chains as well.

We have the following upper bound, due to Halberstam and Richert
  \cite[Theorem 2.4]{HR}.
\begin{lemma} \label{ktuplelemma}
  Given an admissible pattern $\ktuple$ of length $k$,
  the number of integers $x$ such that the $f_i(x)$
  are all simultaneously prime
  and $\max\{f_i(x)\}\le n$ is
  $O(n/(\log n)^k)$.
\end{lemma}
Here the implied constant of the big-$O$ can depend on $k$, but does not
  depend on the coefficients of the $f_i$.

\subsection{Previous Work}

We can obtain a rather straightforward analysis by
  simply finding \textit{all} primes $\le n$,
  and then scanning to see how many tuples match the pattern.
Note that it is not necessary to write all the primes down;
  we can scan them in small batches as we go.
Since the scan phase is fast,
  the bottleneck would be finding the primes using a sieve.
The Atkin-Bernstein sieve \cite{AB2004} does this using at most
  $O(n/\log\log n)$ arithmetic operations and $\sqrt{n}$ space.
Galway \cite{Galway2000} showed how to reduce space further
  to roughly $n^{1/3}$, but this version could not use the wheel
  data structure and requires $O(n)$ arithmetic operations.
See also \cite{Helfgott2020}.

We follow the convention from the analysis of prime sieves
  of not charging for the space used by the output
  of the algorithm.
Note that if we charged for the space of the output,
  we could not expect to do better than $n/(\log n)^{k-1}$
  bits in general.

A more specialized algorithm that searches for only the primes
  in the pattern, and not all the primes, should do much better.
G\"unter L\"oh \cite{Loeh89} and Tony Forbes \cite{Forbes99}
  described algorithms to do exactly this,
  but gave no runtime analysis.
It seems likely their algorithms are faster than $O(n/\log\log n)$,
  but without an analysis, we don't know if this is true or not.
Note that Forbes outlined an odometer-like data structure that
  seems to be similar to the wheel data structure we employ.

Of course, by the prime number theorem, 
  there are only about $n/\log n$ primes $\le n$,
  so the current best prime sieves use $\log n/\log\log n$ arithmetic
  operations per prime on average.
We do not know if anything smaller is possible. 
Applying this average cost to the results of Lemma \ref{ktuplelemma},
  we can hope for an algorithm that takes $O(n/(\log\log n)^k)$
  arithmetic operations to find all primes $\le n$ 
  in a fixed pattern of length $k$.

In essence, this is what we prove as our main result.

\subsection{New Results}

Our contribution is the following.
\begin{theorem}\label{mainresult}
There is an algorithm that,
  when given a list of $k$ distinct linear polynomials over the integers,
  with positive leading coefficients,
  $P=\ktuple$ (the pattern),
  and a search bound $n$,
  finds all integers $x$ such that
  $\max{ \{f_i(x) \} } \le n$ and all the $f_i(x)$ are prime.
This algorithm uses at most
  $O_P(n/(\log\log n)^k)$ arithmetic operations
  and $O(k\sqrt{n})$ bits of space.
\end{theorem}
This algorithm extends the Atkin-Bernstein prime sieve with
  our space-saving wheel sieve.
  See \cite{Sorenson06,Sorenson10a,SW17}.

The $\sqrt{n}$ space needed by this algorithm limits its practicality.
By replacing the Atkin-Bernstein sieve with the sieve of
  Eratosthenes combined with prime tests, we can greatly reduce the need
  for space.

\begin{theorem}\label{spacethm}
Let $c>2$.
There is an algorithm that,
  when given a list of $k>2$ distinct linear polynomials over the integers,
  with positive leading coefficients,
  $P=\ktuple$ (the pattern),
  and a search bound $n$,
  finds all integers $x$ such that
  $\max{ \{f_i(x) \} } \le n$ and all the $f_i(x)$ are prime.
This algorithm uses at most
  $O_P(n/(\log\log n)^{k-1})$ arithmetic operations
  and $n^{1/c}$ bits of space,
  under the assumption of 
  Conjectures \ref{conj1} and \ref{conj2} (see below).
Correctness of the output is unconditional.
\end{theorem}
If $k>6$ then Conjecture \ref{conj1} is not needed.
We use the sieve of Eratosthenes with primes up to
   a bound $B=n^{1/c}$,
  after which we apply a base-2 pseudoprime test and then
  a version of the AKS prime test \cite{AKS04} that 
  was improved to take
  $(\log n)^{6+o(1)}$ time \cite{Lenstra2002}.
Our goal here is to balance the cost of sieving with the cost of prime testing.
For the range $2<k\le 6$,
  to keep the cost of prime testing low enough, 
  we replace the AKS prime test with the
  the pseudosquares prime test of Lukes, Patterson, and Williams \cite{LPW96}.
This prime test takes only $O((\log n)^2)$ time under the condition of
  an unproven but reasonable conjecture 
  on the distribution of pseudosquares due to Bach and Huelsbergen \cite{BH93}.
Note that the correctness of our algorithm 
  does not rely on any unproven conjectures.

The \textit{pseudosquare} $L_p$ is the smallest positive integer that
  is not a square, is $1\bmod 8$, and is a quadratic residue
  modulo every odd prime $q\le p$.
\begin{conjecture}[Bach and Huelsbergen \cite{BH93}] \label{conj1}
  Let $L_p$ be the largest pseudosquare $\le n$.
  Then $p=O(\log n \log\log n)$.
\end{conjecture}

Our second conjecture is needed to bound the cost of
  prime testing after sieving out by primes $\le y$
  in an arithmetic progression.
Let $p(n)$ denote the smallest prime divisor of the integer $n$.
\begin{conjecture}\label{conj2}
  Let $a,b$ be positive integers with $\gcd(a,b)=1$.  Then
$$
  \#\{ n\le x, n\equiv a\bmod{b}, p(n)>y \}
  \ll
  \frac{x}{b}\prod_{p\le y \atop \gcd(p,b)=1} \left(1-\frac{1}{p}\right).
$$
\end{conjecture}
This would be a theorem if $b\le \sqrt{x}$;
  see \cite[(1.7),(1.8)]{Xuan2000,HR}.

We performed a few computations with this second version of our algorithm
  to show its practicality.
A couple of these computational results are new.

The rest of this paper is organized as follows.
In \S\ref{sec:theory} we present our proof of Theorem \ref{mainresult},
  including our model of computation in \S\ref{sec:model},
  a description of our first algorithm in \S\ref{sec:algorithm},
  and its running time analysis in \S\ref{sec:analysis}.
In \S\ref{sec:practice} we discuss our
  second algorithm in \S\ref{sec:algorithm2},
  and its analysis in \S\ref{sec:analysis2},
  thereby proving Theorem \ref{spacethm}. 
We present our computational results in \S\ref{sec:computations},
  including our work on twin primes in \S\ref{sec:twin},
  our work on prime quadruplets in \S\ref{sec:quad},
  and our results on Cunningham chains in \S\ref{sec:chains}.
We wrap up with a discussion of possible future work
  in \S\ref{sec:future}.

\section{Theory\label{sec:theory}}

\subsection{Model of Computation\label{sec:model}}

Our model of computation is a standard random access machine with
  infinite, direct-access memory.
Memory can be addressed at the bit level or at the word level,
  and the word size is $\Theta(\log n)$ bits, if $n$ is the input.
Arithmetic operations on integers of $O(\log n)$ bits take constant time,
  as do memory/array accesses, comparisons, and other basic operations.

We count space used in bits, 
  and we do not include the size of the output.

\subsection{Our First Algorithm\label{sec:algorithm}}

In this section we present the version of our algorithm
  with the smallest running time;
  we perform the analysis in the next section.

The input to the algorithm is 
  the search bound $n$ and
  the pattern,
  which consists of the value of $k$ and
  the list of linear polynomials
  $\ktuple$.
We write $a_i$ for the multiplier and $b_i$ for the offset for
  each form $f_i$.
For simplicity, we often assume that $a_1=1$ and $b_1=0$,
  but this convenience is not required to obtain our complexity bound.
So for example, for Cunningham chains of the first kind we would have
  $a_1=1, a_2=2, a_3=4, \ldots, a_k=2^{k-1}$ and
  $b_1=0, b_2=1, b_3=3, \ldots, b_k=a_k-1$.

\begin{enumerate}
\item
  We begin by finding the list of primes up to $\floor{\sqrt{n}}$
  and choosing a subset to be the \textit{wheel primes} $\Wprimes$.
  Define $W:=\prod_{p\in{\Wprimes}} p$.
  Generally, we put all primes up to a bound $y$
    into $\Wprimes$, with $y$ maximal such that 
    that $W \le \sqrt{n}$.
  Then by the prime number theorem \cite{HW} we have
  $$
   (1/2)\log n \sim \log W  = \sum_{p\le y} \log p \sim y.
  $$
  This implies $\sqrt{n}/\log n \ll W \le \sqrt{n}$.


  In practice, if there is a prime $\le y$ that provides poor filtering
    for the pattern, we consider dropping it from $\Wprimes$ and 
    increasing $y$ to include another prime.

  We must also check all primes $\le y$ to see if they participate
    in the prime pattern.

\item
  Next, we construct the wheel data structure so that it will generate
    all acceptable residues modulo $W$.

  The data structure is initialized with a list of pairwise coprime moduli,
    and for each such modulus $m$, a list of acceptable residues mod $m$,
    encoded as a bit vector of length $m$.
  The wheel modulus $W$ is then the product of these moduli.
  Once initialized,
    the wheel has the ability to enumerate suitable residues modulo $W$
    in amortized constant time per residue.
  The residues do not appear in any particular order.

  Therefore,
    for each prime $p\in{\Wprimes}$
    we compute a bit vector (\texttt{ones[]}) 
    that encodes the list of acceptable residues.
  For any integral $x$, we want $f_i(x)=a_ix+b_i$ to not be divisible
    by $p$.
  So if $p$ divides $a_ix+b_i$, or equivalently if
    $p$ does not divide $a_i$ and so
    $x\equiv -b_i \cdot a_i^{-1} \bmod p$, then bit position $x$ must be
    a zero.
  \begin{quote}\tt
    vector<bool> ones($p,1$); // length $p$, all $1$s to start \\
    for($i=0$; $i<k$; $i=i+1$) \\
    \hspace*{2em} if($a_i \bmod p \ne 0$) \\
    \hspace*{2em} \hspace*{2em} ones[ (-$b_i \cdot a_i^{-1} \bmod p$) ]$=0$;
  \end{quote}
  Continuing the example above, for Cunningham chains of the first kind,
    $p=3$ gives the vector \texttt{001}, and $p=7$ gives the
    vector \texttt{0010111}.

  We then construct the wheel data structure as described in
    \cite[\S4]{Sorenson06}.

\item
  For each residue $r\bmod W$ generated by the wheel,
  we sieve $k$ arithmetic progressions for primes up to $n$,
  $f_i(r)=a_ir+b_i \bmod W$, or $(a_ir+b_i)+j\cdot a_iW$
  for $j=0,\ldots, \lfloor n/(a_iW) \rfloor$.
  We do this using the Atkin-Bernstein sieve.
  (See \cite[\S5]{AB2004} for how to use the sieve to find
  primes in an arithmetic progression.)

  For each $r$, this yields $k$ bit vectors of length $\le n/W$ which are
  combined using a bitwise AND operation to obtain the
  bit positions for where the pattern is satisfied by primes.

\end{enumerate}

\subsubsection*{Example}
Let us find the prime quadruplets ($k=4$, $P=(x,x+2,x+6,x+8)$) 
  up to $n=1000$ using a wheel of
  size $2\cdot3\cdot5\cdot7=210$.

We generate just the one quadruplet $5,7,11,13$ that contains wheel primes.

The only acceptable residue mod $2$ is $1$, for $3$ it is $2$
  (if $x\equiv 1\bmod 3$ then $x+2$ is divisible by $3$),
  giving $5\bmod 6$,
  and for $5$ it is $1$, giving $11\bmod 30$.
For $7$ there are three acceptable residues, $2,3,4$,
  giving us $11, 101, 191 \bmod 210$ by the Chinese remainder theorem.

For $r=11\bmod 210$, we sieve the progression $11+j\cdot 210$ for primes,
  getting the bit vector 10110;
  11, 431, and 641 are prime, 221 and 851 are not.
We then sieve $r+2=13\bmod 210$ for primes, getting 11111; they are all prime.
Sieving $r+6=17\bmod 210$ gives 11011, and
  $r+8\bmod 210$ gives 11101.
$$
\begin{array}{ccccc}
$11$&\sout{221}&$431$&641&\sout{851} \\
$13$&$223$&$433$&$643$&$853$ \\
$17$&$227$&\sout{437}&$647$&$857$ \\
$19$&$229$&$439$&\sout{649}&$859$
\end{array}
$$
Bitwise AND-ing these four vectors together gives 10000.
The only quadruplet we find is $(11,13,17,19)$.

Doing the same for $r=101$, we get the following,
$$
\begin{array}{ccccc}
101&311&521&\sout{731}&941 \\
103&313&523&733&\sout{943} \\
107&317&\sout{527}&\sout{737}&947 \\
109&\sout{319}&\sout{529}&739&\sout{949}
\end{array}
$$
which gives us the quadruplet $(101,103,107,109)$.

Finally for $r=191$ we get
$$
\begin{array}{cccc}
191&401&\sout{611}&821 \\
193&\sout{403}&613&823 \\
197&\sout{407}&617&827 \\
199&409&619&829
\end{array}
$$
and we get two quadruplets, $(191,193,197,199)$ and $(821,823,827,829)$.

%
%
%
%
%

\subsection{Complexity Analysis\label{sec:analysis}}

We'll look at the cost for each step of the algorithm above.

\begin{enumerate}
\item
  We can use the Atkin-Bernstein sieve to find the primes
    up to $\sqrt{n}$ in $O(\sqrt{n}/\log\log n)$ arithmetic operations
    using $n^{1/4}$ space.

\item
  Recall that the largest wheel prime is roughly $(1/2)\log n$.
  Constructing the bit vector \texttt{ones[]} 
    for one prime takes $O(\log n)$ time to initially write all ones,
    and then $O(k)$ time to mark the zeros.
  Summing over all wheel primes gives $O((\log n)^2/\log\log n)$ operations.

  From \cite[Theorem 4.1]{Sorenson06} the total cost to build the
    wheel is $O((\log n)^3)$ operations and it occupies
    $O((\log n)^3/\log\log n)$ space.

\item
  The Atkin-Bernstein sieve finds all primes in an arithmetic progression
    in an interval of size $\sqrt{n}$ or larger in time linear 
    in the length of the interval
    using space proportional to $\sqrt{n}$ \cite[\S5]{AB2004}.
  Therefore,
    sieving for primes takes $O(n/W)$ operations for each of the
    $k$ residue classes $f_i(r) \bmod W$,
    for a total of $O(kn/W)$.
  The cost to generate each value of $r$ using the wheel is
    negligible in comparison.
  The space used is $O(k\sqrt{n})$ bits.

  Next we show that the total number of residues is
    roughly asymptotic to $W/(\log\log n)^k$.
  For a pattern $P$ of size $k$, all but finitely many primes $p$ will have
    $p-k$ possible residues.
  Let $b(P)$ be a constant, depending on the pattern $P$, such that
    all primes larger than $b(P)$ have $p-k$ residues.
  Recall that $y$ is a bound for the largest prime in the wheel.
  Then the total number of residues $r \bmod W$ will be asymptotically
    bounded by
  \begin{eqnarray*}
    \prod_{p\le b(P)} p \cdot \prod_{b(P)<p\le y} (p-k) 
      & \le &  W 
        \prod_{b(P)<p\le y} \frac{p-k }{p} \\
      & = &  W 
        \prod_{b(P)<p\le y} \left(1-\frac{k}{p}\right).
  \end{eqnarray*}
  By Bernoulli's inequality we have $1-k/p \le (1-1/p)^k$, so that
  \begin{eqnarray*}
       W 
        \prod_{b(P)<p\le y} \left(1-\frac{k}{p}\right) 
      & \le & 
       W 
        \prod_{b(P)<p\le y} \left(1-\frac{1}{p}\right)^k \\
     &=&
       W 
        \prod_{p\le b(P)} \left(1-\frac{1}{p}\right)^{-k} \cdot
        \prod_{p\le y} \left(1-\frac{1}{p}\right)^k . \\
  \end{eqnarray*}
  The first product is bounded by a constant that depends on $P$,
    which we refer to as $c(P)$ going forward.
  If $y\ge 4$, then by \cite[(3.26)]{RS62} we know that
    $$ \prod_{p\le y} \left(1-\frac{1}{p}\right) < \frac{1}{\log y}.$$
  As shown above, if $y$ is the largest prime in $\Wprimes$, 
    we have $y\sim (1/2)\log n$.
  Asymptotically, $y\ge (1/10)\log n$ is a very safe underestimate.
  Pulling this all together, we obtain the bound
  $$ c(P) \frac{W}{(\log\log n)^k}  $$
  for the total number of residues $r\bmod W$.
  \nocite{HW}
  
  Multiplying the sieving cost by the number of residues gives \\
  $O(c(P) n/(\log\log n)^k)$ operations.
\end{enumerate}

We have proved Theorem \ref{mainresult}.

It is possible to pin down $c(P)$, 
  the constant that depends on $P$ (and $k$).
We can use \cite[(3.30)]{RS62} and, assuming $b(P)\ge 3$
  (note that $b(P)\ge k$), we have
$$
  \prod_{p\le b(P)} \frac{p}{p-1} < 3.26 \log b(P) .
$$
Bringing in the factor of $k$ introduced in the sieving,
we have
$$
 c(P) \le k ( 3.26 \log b(P) )^{k}.
$$
It remains to get a bound for $b(P)$.
Set $A:=\max\{ |a_j|,|b_j| \}$, an upper bound on all the coefficients
  in the pattern.
Set $F(x):=\prod_{i=1}^k f_i(x) =\prod_{i=1}^k (a_ix+b_i)$.
For a prime $p$, if $F$ has repeated roots modulo $p$, then $p\le b(P)$.
But then $a_ix+b_i\equiv a_jx+b_j \bmod p$ for some $i,j$ with $i\ne j$,
  which means $p$ divides $a_ib_j-a_jb_i$,
  and hence $p\le 2A^2$.
Thus $b(P)\le 2A^2$, and we have the bound
$$
 c(P) \le k ( 3.26 \log (2A^2) )^{k}.
$$
For specific patterns, the constant can be computed explicitly,
  typically giving much better results than this.
In particular, we assumed every prime $\le b(P)$ contributed all its residues
  (which happens if a prime divides all the $a_i$s),
  which is unusual.

As an example, for our pattern $P=(x,x+2,x+6,x+8)$, we have $b(P)=3$,
  and our constant would be
$$
  4 \cdot \left( \frac{1/2}{(1-1/2)^4} \right)
    \left( \frac{1/3}{(1-1/3)^4} \right) = 2\cdot 3^3 = 54.
$$
Here the $1/2$ and $1/3$ multipliers are included because $2$ and $3$
  have only one residue each, instead of $2$ and $3$ respectively.

Note that one of the anonymous referees deserves the primary credit for this
  bound on $c(P)$.

\section{Practice\label{sec:practice}}

The primary difficulty in reaching large values of $n$ with our first
  algorithm is the amount of space it requires.
One way to address this is to create a larger wheel,
  sieve more but shorter arithmetic progressions for primes,
  and rely less on sieving and more on primality tests 
  (in the style of \cite{Sorenson06}) when searching for $k$-tuples.

We use the sieve of Eratosthenes instead of the Atkin-Bernstein sieve
  for the arithmetic progressions, and this is the source of
  the $\log\log n$ factor slowdown.
The gains here are less space needed by a factor of $k$,
  and the effective trial division performs a quantifiable filtering
  out of non-primes.

Instead of sieving by all primes up to $\sqrt{n}$, we sieve only by
  primes up to a bound $B:= n^{1/c}$ for a constant $c>2$.
In practice, we choose $B$ so everything, including space for
  the wheel data structure, the bit vector for sieving,
  and the list of primes $\le B$,
  fits in CPU cache.
We then use a base-2 pseudoprime test,
  followed by a prime test as needed.
For smaller $k$, we use 
  the pseudosquares prime test of 
  Lukes, Patterson, and Williams \cite{LPW96},
  which is fast and deterministic, assuming a sufficient table of
  pseudosquares is available. 
Importantly, it takes advantage of
  the trial division effect of the sieve of Eratosthenes.
For larger $k$, we can simply use the AKS prime test \cite{AKS04}.

This change means we can get by with only $O(B)$ space. 
Choosing $B$ larger or smaller
  gives a tradeoff between the cost of sieving and the cost of
  performing base-2 pseudoprime tests.

\subsection{Our Second Algorithm\label{sec:algorithm2}}

\begin{enumerate}
\item
  Choose a constant $c>2$, and then
    set $B:=2^{\lfloor (1/c)\log_2 n \rfloor}$,
    a power of $2$ near $n^{1/c}$.
  We begin by finding the list of primes up to $B$
  and dividing them into the two sets $\Wprimes$ and $\Sprimes$.
  Small primes go into $\Wprimes$ and the remainder go in $\Sprimes$.
  We want $W:=\prod_{p\in{\Wprimes}} p$ to be
    as large as possible with the constraint that $W \le n/B$.
  This implies that the largest prime in $\Wprimes$ will be
    roughly $\log n$.

\item
  If $k\le 6$, we will need to perform the pseudosquares prime test,
  so in preparation, find all pseudosquares $L_p\le n/B$
  (see \cite{LPW96}).

\item
  Next, as before, 
  we construct the wheel data structure so that it will generate
    all possible correct residues modulo $W$.

\item
  For each residue $r\bmod W$ generated by the wheel,
    we construct a bit vector \texttt{v[]} of length $n/W$.
  Each vector position \texttt{v[$j$]}, for $j=0,\ldots,\lfloor W/n \rfloor$,
    represents $x(j)=r+j\cdot W$
    for the $k$-tuple $(f_1(x(j)), f_2(x(j)), \ldots, f_k(x(j)))$.
  We initialize \texttt{v[$j$]}$=1$,
    but clear it to $0$ if we find a prime $p\in\Sprimes$
    where $p\mid f_i(x(j))$ for some $i$.

  \begin{quote}\tt
  \begin{tabbing}MM\=MM\=MM\=MM\=\kill
    for( $p\in \Sprimes$) \+\\
      $winv := W^{-1} \bmod p$; \\
      for($i:=0$; $i<k$; $i:=i+1$) \+\\
         $j:=winv\cdot(-b_ia_i^{-1}-r) \bmod p$; \\
         while($j<n/W$) \+\\
           \texttt{v[$j$]}$:=0$; \\
           $j:=j+p$; \\
  \end{tabbing}
  \end{quote}
    

  Once this sieving is complete, the only integers $j$
    with \texttt{v[$j$]}$=1$
    that remain, satisfy the property that all the
    $f_i(x(j))$ have no prime divisors less than $B$.

\item
  For each such $x(j)$ remaining (that is, \texttt{v[$j$]}$=1$), 
    we first do a base-2 strong pseudoprime test on $f_1(x(j))$.  
  If it fails, we cross it off (set \texttt{v[$j$]}$=0$).
  If it passes, we try $f_2(x)$ and so forth, 
    keeping \texttt{v[$j$]}$=1$ only if all $k$
    values $f_i(x(j))$ pass the pseudoprime test.
  We then perform a full prime test on the $f_i(x(j))$ for all $i$.
  If $k\le 6$, we use
    the Lukes, Patterson, and Williams
    pseudosquares prime test \cite{LPW96} as done in \cite{Sorenson06}.
  For larger $k$, we use the AKS prime test \cite{AKS04}.
  (This is for the purposes of the theorem; in practice, 
    the pseudosquares prime test is faster, so we use that instead.)
  If all the $f_i(x(j))$ pass the prime tests, 
    the corresponding $k$-tuple is written for output.

\end{enumerate}
This version of the algorithm works best for $k\ge4$.
When $k\le 3$,
  the prime tests become the runtime bottleneck, and so
  we recommend using $B=\sqrt{n}$ so that the base-2 pseudoprime tests and
  the pseudosquares prime test are not needed, as the sieving
  will leave only primes. 

\subsubsection*{Example}
We use the same example as above, finding prime quadruplets up to $5000$,
  which uses the pattern $(x,x+2,x+6,x+8)$.
We go a bit higher this time to illustrate the sieving.
Recall our wheel uses the primes $2\cdot3\cdot5\cdot7$,
  and generates the three residues $11,101,191$ modulo $210$.

We will use $B=20$ as our sieve bound, so $\Sprimes=\{11,13,17,19\}$.
Like with the wheel primes,
  quadruplets that include sieve primes must be generated separately,
  so we output the quadruplet $(11,13,17,19)$ at this point.

For $r=11$, we have our progression
$$ \begin{array}{cccccccccc}
      11 &  221 &  431 &  641 & 851
  & 1061 & 1271 & 1481 & 1691 & 1901 \\
    2111 & 2321 & 2531 & 2741 & 2951 
  & 3161 & 3371 & 3581 & 3791 & 4001 \\
    4211 & 4421 & 4631 & 4841
\end{array} $$
For $p=11$, we cross off integers that are $0,3,5,9$ modulo $11$.
This takes four passes: 
  on the first pass we remove $11,2321,4631$, which are $0\bmod 11$;
  on the second pass we remove $641,2951$, which are $3\bmod 11$
  (for example $649=11\cdot59$);
  the third removes $1061,3371$, which are $5\bmod 11$;
  and the fourth removes $1901,4211$, which are $9\bmod 11$.
Observe that for a given residue modulo $11$, the numbers to cross off
  are exactly $11$ spaces apart.
$$ \begin{array}{cccccccccc}
      \sout{11}&  221 &  431 &  \sout{641}& 851
  & \sout{1061}& 1271 & 1481 & 1691 & \sout{1901} \\
    2111 & \sout{2321}& 2531 & 2741 & \sout{2951}
  & 3161 & \sout{3371}& 3581 & 3791 & 4001 \\
    \sout{4211}& 4421 & \sout{4631}& 4841
\end{array} $$
For $p=13$, we cross off integers that
  are $0,5,7,11$ modulo $13$. 
We remove $221,2951$ for $0\bmod 13$,
$2111,4841$ for $5\bmod 13$,
$2321$ for $7\bmod 13$,
and $11,2741$ for $11\bmod 13$.
$$ \begin{array}{cccccccccc}
      \sout{11}&  \sout{221}&  431 &  \sout{641}& 851
  & \sout{1061}& 1271 & 1481 & 1691 & \sout{1901} \\
    \sout{2111}& \sout{2321}& 2531 & \sout{2741}& \sout{2951}
  & 3161 & \sout{3371}& 3581 & 3791 & 4001 \\
    \sout{4211}& 4421 & \sout{4631}& \sout{4841}
\end{array} $$
For $p=17$, we cross off
  integers that are $0,9,11,15$ modulo $17$.
We remove $221,3791$ for $0\bmod 17$,
$2321$ for $9\bmod 17$,
$2531$ for $15\bmod 17$,
and $11,3581$ for $11\bmod 17$.
$$ \begin{array}{cccccccccc}
      \sout{11}&  \sout{221}&  431 &  \sout{641}& 851
  & \sout{1061}& 1271 & 1481 & 1691 & \sout{1901} \\
    \sout{2111}& \sout{2321}& \sout{2531}& \sout{2741}& \sout{2951}
  & 3161 & \sout{3371}& \sout{3581}& \sout{3791}& 4001 \\
    \sout{4211}& 4421 & \sout{4631}& \sout{4841}
\end{array} $$
For $p=19$, we cross off integers that are $0,11,13,17$ modulo $19$.
We remove $11,4001$ for $11\bmod 19$,
$431,4421$ for $13\bmod 19$,
$1271$ for $17\bmod 19$,
and $1691$ for $0\bmod 19$.
$$ \begin{array}{cccccccccc}
      \sout{11}&  \sout{221}&  \sout{431}&  \sout{641}& 851
  & \sout{1061}& \sout{1271}& 1481 & \sout{1691}& \sout{1901} \\
    \sout{2111}& \sout{2321}& \sout{2531}& \sout{2741}& \sout{2951}
  & 3161 & \sout{3371}& \sout{3581}& \sout{3791}& \sout{4001} \\
    \sout{4211}& \sout{4421}& \sout{4631}& \sout{4841}
\end{array} $$
At this point we perform base-2 strong pseudoprime tests,
  followed by prime tests as needed.
Here 851 and 3161 are not prime, and
 1481 leads to the quadruplet $(1481,1483,1487,1489)$.

This is then repeated with $r=101,191$.

\subsection{Complexity Analysis\label{sec:analysis2}}

Finding the primes up to $B$ takes $O(B)$ time
  using $O(\sqrt{B})$ space, well within our bounds.
See \cite{Sorenson06} for a sublinear time algorithm
  to find all needed pseudosquares.
In practice, all pseudosquares up to $10^{25}$ are known \cite{Sorenson10a}.
The cost in time and space to build the wheel is, up to a constant factor,
  the same.
So we now focus on steps (4) and (5).

As shown above,
  the number of residues to check $\bmod$ $W$ is
  $O_P(W/(\log\log n)^k)$.
The time to sieve each interval of length $n/W$ using primes up to $B$
  is at most proportional to
  $$
  \sum_{p\le B} \frac{kn}{pW} \sim \frac{kn\log\log B}{W} 
    \sim \frac{kn\log\log n}{W} .
  $$
Here the multiplier $k$ is required because we cross off $k$ residues modulo
  most of the primes $p\le B$.
That said, this multiple of $k$ can be absorbed into the implied constant that
  depends on the pattern $P$, $c(P)$, from earlier.

At this point we make use of Conjecture \ref{conj2} to bound
  the number of integers free from prime divisors $\le B$ in an arithmetic
  progression.

With this assumption in hand, by Mertens's theorem,
  an average of at most
  $$
    \frac{n}{W} \prod_{y<p\le B} \left( 1- \frac{1}{p} \right)
    \quad \ll \quad \frac{n\log y}{W \log B}
    \quad \sim \quad \frac{n\log\log n}{W \log n}
  $$
  vector locations remain to be prime tested.
(Note that we cannot make any assumptions about the relative
   independence of the primality of the $f_i(x)$ values for different $i$,
   and so we cannot use a $(1-k/p)$ factor here.)

A single base-2 strong pseudoprime test takes $O(\log n)$
  operations to perform, for a total cost proportional to
  $$ \frac{ kn \log\log n }{ W \log n } \log n
   \sim  \frac{ kn \log\log n }{ W } $$
  arithmetic operations
  to do the base-2 strong pseudoprime tests for each value of $r\bmod W$.
This matches the sieving cost
  of $O_P( n\log\log n / W)$ from above.
(Note that if we deliberately choose a larger value for $B$,
  the increased sieving will decrease the number of pseudoprime tests needed.
 This tradeoff can be used to fine-tune the running time of the algorithm.)

Thus, the total cost for sieving and base-2 pseudoprime tests
  is $$ O_P\left( \frac{n}{(\log\log n)^{k-1}}\right),$$
which we obtain by multiplying by the number of
  residues $O_P(W/(\log\log n)^k)$.

Next we need to count integers that pass the base-2 strong
  pseudoprime test.
Such integers are either prime, or composite base-2 pseudoprimes.
We switch to counting across all residues $r\bmod W$
  to obtain an overall bound.

Lemma \ref{ktuplelemma} tells us that at most
  $O(n/(\log n)^k)$ integers are prime that fit the pattern,
  so this is an upper bound on primes that pass the base-2 pseudoprime test.


Pomerance \cite{Pomerance81} showed that
  the number of composite base-2 pseudoprimes $\le n$ is bounded by 
$$
  {n}{e^{-\sqrt{\frac{\log n \log\log\log n }{ \log\log n}}}}
  \ll \frac{n}{(\log n)^{k+1}}
$$ 
which is negligible.
This plus the bound for primes above gives us the $O(n/(\log n)^k)$ bound we
  desire for all integers that pass the base-2 pseudoprime test.

Next, to bound the cost of prime tests, we have two cases:
  $k>6$, or $2< k \le 6$.

For $k>6$, we use the AKS prime test \cite{AKS04}
  (improved in \cite{Lenstra2002})
  which takes time $O((\log n)^{6+o(1)})$.
The cost of applying the AKS prime test
  to all the integers $f_i(x)$ after they all pass a base-2 pseudoprime test
  is at most proportional to
  $$ k\cdot (\log n)^{6+o(1)} \cdot \frac{n}{(\log n)^k}
   \ll \frac{kn}{(\log n)^{k-6+o(1)}}$$
  which is $o(kn/(\log\log n)^{k})$ for $k>6$.

Note that when $k$ is large, in practice we might only do the
  base-2 pseudoprime tests, and then run full prime tests on the
  output afterwards, since the amount of output will be rather small.
 
For $2<k\le 6$, 
  Conjecture \ref{conj1}
  implies that the pseudosquares prime test
  takes $O((\log n)^2)$ arithmetic operations to test
  integers $\le n$ for primality,
  given a table of pseudosquares $\le n$.
If $n$ has no prime divisors below $B$, then pseudosquares up to $n/B$
  suffice.  See \cite{LPW96,Sorenson06}.

So, under the assumption of Conjecture \ref{conj1},
  the cost of applying the pseudosquares prime test
  to all the integers $f_i(x)$ after they all pass a base-2 pseudoprime test
  is at most proportional to
  $$
  k\cdot (\log n)^2  \cdot  \frac{n}{(\log n)^k} 
    \ll \frac{kn}{ (\log n)^{k-2} }
  $$
  and this is $o(kn/(\log\log n)^k)$ for $k>2$.

The space used is dominated by the length of the sieve intervals and
  the space needed to store the primes in $\Sprimes$,
  which is $O(B)$ bits.

This completes the proof of Theorem \ref{spacethm}. 

\section{Computations\label{sec:computations}}

As mentioned previously,
  we implemented several versions of our second algorithm
  to see what we could compute.
We looked for new computational records that were within reach of
  our university's nice but aging hardware.
Below we discuss some of the results of those computations.
Some of the implementation details are specific to a particular
  computation.
Here are four remarks about implementation details
  that these computations had in common.

\begin{enumerate}
\item
  We wrote our programs in C++ using the GNU compiler
    under Linux.
  GMP was used for multiprecision arithmetic when necessary.
  Note that it is fairly easy to write the code such that GMP
    was needed only on occasion and for prime tests.
\item
  We used MPI 
    and ran our code on Butler University's cluster \textit{Big Dawg}.
  This machine has 16 compute nodes
    with 12 cores (2 CPUS) each at optimal capacity;
    our average utilization was around 150 of the 192 cores 
    due to compute nodes going down from time to time.
  The CPU is the Intel Xeon CPU E5-2630 0 @ 2.30GHz
    with 15 MB cache, with 6 cores per CPU.

  To parallelize the algorithm, we striped on the residues $r\bmod W$.
  In other words,
    all $\nu$ processors stepped through all the $r \bmod W$ residues,
    but only sieved every $\nu$th residue for primes.
  This meant there was very little communication overhead
    except for when periodic checkpoints were done, 
    about every 15-30 minutes.
\item
  We usually chose our wheel size ($W$) and sieve intervals so that
    the size of each interval ($n/W \approx B$) was at most a few megabytes
    so that it would fit in the CPU cache.
  We used a \texttt{vector<bool>}, which packs bits.
\item
  For larger values of $k$,
    we observed that when sieving by smaller primes $p$ by each of the $k$
    residues,
    we might find that almost all the bits of the current interval were cleared
    long before we reached the sieving limit $B$,
    so we created a simple early-abort strategy that was able to save time.

  The very few remaining bits were tested with the base-2 strong pseudoprime
    test even though we had not sieved all the way to $B$.
  We also, then, replaced the use of the pseudosquares prime test with
    strong pseudoprime tests \cite{Miller76}
    using results from \cite{SW17} so that only a few bases were needed,
    due to the spotty trial-division information.

\item
  We found that, especially for larger $k$, 
    our algorithm spent more time on sieving
    than prime testing.
  As mentioned previously, for $k=3$ the prime testing dominates the
    running time in practice, and it is worthwhile to use $B=\sqrt{n}$
    so that prime testing is not required.
  
\end{enumerate}

\subsection{Twin Primes and Brun's Constant\label{sec:twin}}

Let $\pi_2(X)$ count the twin prime pairs $(p, p+2)$ with $p < X$ 
  and $S_2(X)$ be the sum of the reciprocals of their elements.  
Thomas Nicely computed these functions up to $2\cdot10^{16}$
  (See \texttt{http://www.trnicely.net/\#PI2X}).  
We verified his computations to 14 digits 
  and extended the results to $X = 10^{17}$.  
A portion of our computational results are in the table below.
\begin{center}
\begin{tabular}{ |r |c |c| }
\hline
 $X$ & $\pi_2(x)$ & $S_2(X)$ \\ \hline \hline
$1\cdot10^{16}$ & 10304195697298	& 1.8304844246583 \\  
$2\cdot10^{16}$ & 19831847025792	& 1.8318080634324 \\ \hline
$3\cdot10^{16}$ & 29096690339843	 &1.8325599218628 \\
$4\cdot10^{16}$ & 38196843833352	 &1.8330837014776 \\
$5\cdot10^{16}$ & 47177404870103	& 1.8334845790134 \\
$6\cdot10^{16}$ & 56064358236032	& 1.8338086822020 \\
$7\cdot10^{16}$ & 64874581322443	& 1.8340803303554 \\
$8\cdot10^{16}$ & 73619911145552 &	  1.8343139034256 \\
$9\cdot10^{16}$ & 82309090712061	 &1.8345186031523 \\
$10\cdot10^{16}$ & 90948839353159&	  1.8347006694414 \\ \hline
\end{tabular}
\end{center}
The last section of Klyve's PhD Thesis \cite{Klyve2007}
  describes how to use this information to derive bounds for
  Brun's constant.

We have four remarks on our algorithm implementation:
\begin{enumerate}
  \item As mentioned above, for small $k$ like $k=2$,
    it is more efficient to set $B=\sqrt{n}$ so that sieving 
    also determines primality, thereby avoiding base-2 strong pseudoprime
    tests and primality tests.
  \item We computed $S_2$ using Kahan summation \cite{Kahan65} with
    the \texttt{long double} data type in C++, which gave us 17 digits,
    14 of which were accurate;
    Thomas Nicely has data with 53 digits of accuracy.
    The partial sums were accumulated in 10,000 buckets
    for each process, and then the buckets were in turn added up
    across processes using Kahan summation.
  \item
    Our computation took roughly 3 weeks of wall time, which included
    at least one restart from a checkpoint.
    Our verification of Nicely's work to $10^{16}$ took 42 hours.
  \item
    We used a wheel with 
      $W=6469693230 
   = 2\cdot 3\cdot 5\cdot 7\cdot 11\cdot 13\cdot 17\cdot 19\cdot 23\cdot 29$.
    Note that this is roughly $20\cdot\sqrt{10^{17}}$.
    There were 
      $214708725$ residues $r\bmod W$ to sieve.
\end{enumerate}
\nocite{Nicely96,Klyve2007}

See \texttt{OEIS.org} sequence A007508,
  the number of twin prime pairs below $10^n$.

\subsection{Quadruplet Primes\label{sec:quad}}

A related sum involves the reciprocals of
 the elements of the prime tuple $(p, p+2, p+6, p+8)$.
Let $\pi_4(X)$ count these tuplets up to $X$,
and let $S_4(X)$ be the sum of the reciprocals of their elements.  
Thomas Nicely computed these functions up to $2\cdot10^{16}$.
We extended this computation and partial results are in the table below.
The first two lines are Thomas Nicely's own results,
  which we verified.
\begin{center}
\begin{tabular}{ |r |r |c| }
\hline
 $X$ & $\pi_4(x)$ & $S_4(X)$ \\ \hline \hline
$1\cdot10^{16}$ & 25379433651 & 0.8704776912340 \\ 
$2\cdot10^{16}$ & 46998268431&	0.8704837109481 \\ \hline
$3\cdot10^{16}$ &67439513530   &0.8704870310432 \\
$4\cdot10^{16}$ &87160212807&	0.8704893020026 \\
$5\cdot10^{16}$ & 106365371168&	0.8704910169467 \\
$6\cdot10^{16}$ & 125172360474 &0.8704923889088 \\
$7\cdot10^{16}$ & 143655957845&	0.8704935288452 \\
$8\cdot10^{16}$ & 161868188061 &0.8704945017556 \\
$9\cdot10^{16}$ & 179847459283&	0.8704953489172 \\
$10\cdot10^{16}$ & 197622677481&0.8704960981105 \\ \hline
\end{tabular}
\end{center}
This computation took about 4 days, and we
  used a separate program rather than looking for
  pairs of twin primes in the first program.
Even though $k=4$ is large enough to use prime tests,
  we found that sieving to $\sqrt{n}$ was faster in practice.

We used a wheel with $W=200560490130
= 2\cdot 3\cdot 5\cdot 7\cdot 11\cdot 13\cdot 17\cdot 19\cdot 23\cdot 29\cdot 31$,
  which gave $472665375$ residues.

See \texttt{OEIS.org} sequence A050258,
  the number of prime quadruplets with largest member $< 10^n$.

\subsection{Cunningham Chains\label{sec:chains}}

We have two new computational results for Cunningham chains.
\begin{enumerate}
\item
We found the smallest chain of length 15 of the first kind,
  and it begins with the prime
$$ p=90616\ 21195\ 84658\ 42219. $$
The next four chains of this length of the first kind 
  begin with
\begin{quote}
$1\ 13220\ 80067\ 50697\ 84839$ \\
$1\ 13710\ 75635\ 40868\ 11919$ \\
$1\ 23068\ 71734\ 48294\ 53339$ \\
$1\ 40044\ 19781\ 72085\ 69169$
\end{quote}
This computation took roughly a month of wall time.
Here we used wheel size $W=19835154277048110
= 2\cdot 3\cdot 5\cdot 7\cdot 11\cdot 13\cdot 17\cdot 19
   \cdot 23\cdot 29\cdot 37\cdot 41\cdot 43\cdot 47$,
  with $12841500672$ residues to sieve.
Note that $31$ is a rather badly behaved prime for larger Cunningham chains
  (only $5$ residues are excluded),
  so we left it out of the wheel.

See \texttt{OEIS.org} sequence A005602,
  smallest prime beginning a Cunningham chain of length $n$
  (of the first kind).

\item
In 2008 Jaroslaw Wroblewski found a Cunningham chain of 
  length 17 of the first kind, starting with 
\[ p = 27\ 59832\ 93417\ 13865\ 93519, \]
and we were able to show that this is in fact the smallest such
  chain of that length.

This computation took roughly three months of wall time.
We used $W=1051263176683549830
= 2\cdot 3\cdot 5\cdot 7\cdot 11\cdot 13\cdot 17\cdot 19
   \cdot 23\cdot 29\cdot 37\cdot 41\cdot 43\cdot 47\cdot 53$,
  with $35864945424$ residues to sieve.
With roughly three times as many residues as the previous computation,
  it took roughly three times as long to complete.

\end{enumerate}

\section{Discussion and Future Work\label{sec:future}}

In summary, we have described and analyzed two algorithms for
  finding primes in patterns, and then shown that
  the second of these algorithms is quite practical by performing a few
  computations.

We have some ideas for future work.
\begin{enumerate}
\item In the Introduction, 
  we mentioned that our algorithms could be used to find
  Carmichael numbers by finding prime triplets that satisfy
  the pattern $(6x+1,12x+1,18x+1)$, but we have not yet done
  that computation \cite{Chernick39}.
\item 
  Does it make sense to use Bernstein's doubly-focused enumeration to
  attempt to further reduce the running time?
  See \cite{Bernstein04,Sorenson10a,WW2006}
\item 
  A natural extension to our algorithms here is to allow the
  linear polynomials $f_i$ to potentially be higher degree, 
  irreducible polynomials.
  See Schinzel's Hypothesis H (See \cite{SS58} and \cite[\S1.2.2]{CP})
    and the Bateman-Horn conjecture \cite{BH62}.
\item
  An algorithm for twin primes with space roughly $n^{1/3}$
    that runs in $O(n/(\log\log n)^2)$ time would be nice.
\end{enumerate}
\nocite{UPINT}

\section*{Acknowledgements}

We wish to thank Frank Levinson for his generous gift that funds the
  Butler cluster supercomputer \textit{Big Dawg},
and the Holcomb Awards Committee for their financial support of this work.

Thanks to Thom\`{a}s Oliveira e Silva for several helpful comments
  on an earlier version of this paper.

We also wish to thank two anonymous referees
  for many detailed, helpful comments
  on earlier versions of this paper.
In particular, one of the referees contributed the analysis of
  the implied constant in \S2.3.

A preliminary version of this work was presented as a poster
  at the ANTS XIII conference in Madison, Wisconsin in July of 2018.


\begin{thebibliography}{10}

\bibitem{AKS04}
Manindra Agrawal, Neeraj Kayal, and Nitin Saxena.
\newblock P{RIMES} is in {P}.
\newblock {\em Ann. of Math. (2)}, 160(2):781--793, 2004.

\bibitem{AB2004}
A.~O.~L. Atkin and D.~J. Bernstein.
\newblock Prime sieves using binary quadratic forms.
\newblock {\em Mathematics of Computation}, 73:1023--1030, 2004.

\bibitem{BH93}
Eric Bach and Lorenz Huelsbergen.
\newblock Statistical evidence for small generating sets.
\newblock {\em Math. Comp.}, 61(203):69--82, 1993.

\bibitem{BH62}
Paul~T. Bateman and Roger~A. Horn.
\newblock A heuristic asymptotic formula concerning the distribution of prime
  numbers.
\newblock {\em Math. Comp.}, 16:363--367, 1962.

\bibitem{Bernstein04}
Daniel~J. Bernstein.
\newblock Doubly focused enumeration of locally square polynomial values.
\newblock In {\em High primes and misdemeanours: lectures in honour of the 60th
  birthday of Hugh Cowie Williams}, volume~41 of {\em Fields Inst. Commun.},
  pages 69--76. Amer. Math. Soc., Providence, RI, 2004.

\bibitem{Chernick39}
Jack Chernick.
\newblock On {F}ermat's simple theorem.
\newblock {\em Bull. Amer. Math. Soc.}, 45(4):269--274, 1939.

\bibitem{CP}
R.~Crandall and C.~Pomerance.
\newblock {\em Prime Numbers, a Computational Perspective}.
\newblock Springer, 2001.

\bibitem{Dickson04}
L.~E. Dickson.
\newblock A new extension of {Dirichlet's} theorem on prime numbers.
\newblock {\em Messenger of mathematics}, 33:155--–161, 1904.

\bibitem{Forbes99}
Tony Forbes.
\newblock Prime clusters and {C}unningham chains.
\newblock {\em Math. Comp.}, 68(228):1739--1747, 1999.

\bibitem{Galway2000}
William~F. Galway.
\newblock Dissecting a sieve to cut its need for space.
\newblock In {\em Algorithmic number theory (Leiden, 2000)}, volume 1838 of
  {\em Lecture Notes in Comput. Sci.}, pages 297--312. Springer, Berlin, 2000.

\bibitem{GR98}
Daniel~M. Gordon and Gene Rodemich.
\newblock Dense admissible sets.
\newblock In {\em Algorithmic number theory ({P}ortland, {OR}, 1998)}, volume
  1423 of {\em Lecture Notes in Comput. Sci.}, pages 216--225. Springer,
  Berlin, 1998.

\bibitem{UPINT}
Richard~K. Guy.
\newblock {\em Unsolved problems in number theory}.
\newblock Problem Books in Mathematics. Springer-Verlag, New York, third
  edition, 2004.

\bibitem{HR}
H.~Halberstam and H.-E. Richert.
\newblock {\em Sieve Methods}.
\newblock Academic Press, 1974.

\bibitem{HL23}
G.~H. Hardy and J.~E. Littlewood.
\newblock Some problems of `partitio numerorum'; iii: On the expression of a
  number as a sum of primes.
\newblock {\em Acta Mathematica}, 44(1):1--70, Dec 1923.

\bibitem{HW}
G.~H. Hardy and E.~M. Wright.
\newblock {\em An Introduction to the Theory of Numbers}.
\newblock Oxford University Press, 5th edition, 1979.

\bibitem{Helfgott2020}
Harald~Andr\'{e}s Helfgott.
\newblock {\em Mathematics of Computation}, 89:333--350, 2020.

\bibitem{Kahan65}
W.~Kahan.
\newblock Pracniques: Further remarks on reducing truncation errors.
\newblock {\em Commun. ACM}, 8(1):40--, January 1965.

\bibitem{Klyve2007}
Dominic Klyve.
\newblock {\em Explicit Bounds on Twin Primes and {Brun's} Constant}.
\newblock {PhD} in mathematics, Dartmouth College, Hanover, NH USA, 2007.

\bibitem{Lenstra2002}
H.~W. Lenstra, Jr. and Carl Pomerance.
\newblock Primality testing with "gaussian" periods.
\newblock In M.~Agrawal and A.~Seth, editors, {\em Proceedings of the 22nd
  Conference on Foundations of Software Technology and Theoretical Computer
  Science}, FST TCS '02, Berlin, Heidelberg, 2002. Springer-Verlag.
\newblock LNCS 2556.

\bibitem{Loeh89}
G\"unter L\"oh.
\newblock Long chains of nearly doubled primes.
\newblock {\em Math. Comp.}, 53(188):751--759, 1989.

\bibitem{LPW96}
R.~F. Lukes, C.~D. Patterson, and H.~C. Williams.
\newblock Some results on pseudosquares.
\newblock {\em Math. Comp.}, 65(213):361--372, S25--S27, 1996.

\bibitem{Miller76}
G.~Miller.
\newblock Riemann's hypothesis and tests for primality.
\newblock {\em Journal of Computer and System Sciences}, 13:300--317, 1976.

\bibitem{Nicely96}
T.~Nicely.
\newblock Enumeration to $10^{14}$ of the twin primes and {Brun}'s constant.
\newblock {\em Virginia J. Sci.}, 46:195--204, 1996.

\bibitem{Pomerance81}
Carl Pomerance.
\newblock On the distribution of pseudoprimes.
\newblock {\em Math. Comp.}, 37(156):587--593, 1981.

\bibitem{Riesel}
Hans Riesel.
\newblock {\em Prime Numbers and Computer Methods for Factorization}, volume
  126 of {\em Progress in Mathematics}.
\newblock Birkh{\"a}user, 2nd edition, 1994.

\bibitem{RS62}
J.~B. Rosser and L.~Schoenfeld.
\newblock Approximate formulas for some functions of prime numbers.
\newblock {\em Illinois Journal of Mathematics}, 6:64--94, 1962.

\bibitem{SS58}
A.~Schinzel and W.~Sierpi\'nski.
\newblock Sur certaines hypoth\`eses concernant les nombres premiers.
\newblock {\em Acta Arith.}, 4:185--208, 1958.
\newblock Erratum 5 (1958), 259.

\bibitem{Sorenson06}
Jonathan~P. Sorenson.
\newblock The pseudosquares prime sieve.
\newblock In Florian Hess, Sebastian Pauli, and Michael Pohst, editors, {\em
  Proceedings of the 7th International Symposium on Algorithmic Number Theory
  (ANTS-VII)}, pages 193--207, Berlin, Germany, July 2006. Springer.
\newblock LNCS 4076, ISBN 3-540-36075-1.

\bibitem{Sorenson10a}
Jonathan~P. Sorenson.
\newblock Sieving for pseudosquares and pseudocubes in parallel using
  doubly-focused enumeration and wheel datastructures.
\newblock In Guillaume Hanrot, Francois Morain, and Emmanuel Thom\'e, editors,
  {\em Proceedings of the 9th International Symposium on Algorithmic Number
  Theory (ANTS-IX)}, pages 331--339, Nancy, France, July 2010. Springer.
\newblock LNCS 6197, ISBN 978-3-642-14517-9.

\bibitem{SW17}
Jonathan~P. Sorenson and Jonathan Webster.
\newblock Strong pseudoprimes to twelve prime bases.
\newblock {\em Math. Comp.}, 86(304):985--1003, 2017.

\bibitem{WW2006}
Kjell Wooding and H.~C. Williams.
\newblock Doubly-focused enumeration of pseudosquares and pseudocubes.
\newblock In {\em Proceedings of the 7th International Algorithmic Number
  Theory Symposium (ANTS VII)}, Berlin, Germany, 2006.

\bibitem{Xuan2000}
T.~Z. Xuan.
\newblock Integers free of small prime factors in arithmetic progressions.
\newblock {\em Nagoya Mathematical Journal}, 157:103--127, 2000.

\bibitem{Zhang14}
Yitang Zhang.
\newblock Bounded gaps between primes.
\newblock {\em Ann. of Math. (2)}, 179(3):1121--1174, 2014.

\end{thebibliography}

\end{document}